\documentclass{amsart}

\newcommand{\<}{\langle}
\renewcommand{\>}{\rangle}

\renewcommand{\d}{\delta}

\newcommand{\Z}{\mathbb Z}

\newcommand{\C}{\mathbb C}

\newcommand{\G}{ G}

\newcommand{\cstar}{\C^\times}
\newcommand{\Cstar}{\C^\times}

\theoremstyle{plain}

\theoremstyle{definition}

\theoremstyle{remark}

\begin{document}

\title[Central Extension]{Yet another construction of the 
central extension of the loop group. }
\author{Michael K. Murray}
\address[Michael K. Murray]
{Department of Pure Mathematics\\
University of Adelaide\\
Adelaide, SA 5005 \\
Australia}
\email[Michael K. Murray]{mmurray@maths.adelaide.edu.au}
\thanks{The first author acknowledges the support of the Australian
Research Council.}
\author{Daniel Stevenson}
\email[Daniel Stevenson]{dstevens@maths.adelaide.edu.au}
\thanks{The second author acknowledges the support of the Australian
Research Council.}


\begin{abstract} We give  a characterisation of central extensions of 
a Lie group $G$ by $\Cstar$ in terms of a differential two-form
on $G$ and a differential one-form on $G \times G$.  This is applied
to the case of the central extension of the loop group.
\end{abstract}
\maketitle
\section{Introduction} 
Let $G$ and $A$ be groups.  A {\em central extension} of 
$G$ by $A$ is another group $\hat G$ and a homomorphism
$\pi\colon\hat G \to G$ whose kernel is isomorphic to $A$ and in the 
center of $\hat G$.  Note that because $A$ is in the center
of $\hat G$ it is necessarily abelian. We will be interested ultimately
in the case that $G = \Omega(K)$ the {\em loop group} of all
smooth maps from the circle $S^1$ to a Lie group $K$ with
pointwise multiplication but the theory developed applies
to any Lie group $G$.  

\section{Central extension of groups}
Consider first the case when $G$ is just a group
and  ignore questions of continuity or differentiability. 
In this case we can choose  a {\em section} of the map $\pi$. That
is a map $s \colon G \to \hat G$ such that $\pi(s(g)) = g $
for every $g \in G$.  From this section we can construct a bijection
$$
\phi \colon A \times G \to \hat G
$$
by $\phi(g, a) = \iota(a)s(g)$ where $\iota \colon A \to \hat G$
is the identification of $A$ with the kernel of $\pi$.  So we know 
that  as a set $\hat G$ is just the product $A \times G$. However as a 
group $\hat G$ is not generally a product. To see what it is note that
$\pi(s(g)s(h)) =\pi(s(g))\pi(s(h)) = gh = \pi(s(gh))$ so that 
$s(g)s(h) = c(g, h) s(gh)$ where $c \colon G \times G \to A$. The 
bijection $\phi \colon A \times G \to \hat G$ induces a 
product on $A \times G$ for which $\phi$ is a homomorphism. 
To calculate this  product we note that
 \begin{align*}
     \phi(a, g) \phi(b, h) &= \iota(a)s(g)\iota(b)s(h) \\
                          &= \iota(ab)s(g)s(h)\\
			  &= \iota(ab)c(gh)s(gh).
			  \end{align*}
Hence the product on $A \times G$ is given by 
$(a, g) \star (b, h) = (abc(g,h) gh)$ and the map $\phi$ is a 
group isomorphism between $\hat G$ and $A \times G$ with this
product. 

Notice that if we choose a different section $\tilde s$ then 
$\tilde s = sh$ were $h \colon G \to A$. 

It is straightforward to check that if we pick any $c \colon G \times 
G \to A$ and define a product on $A \times G$ by
$(a, g) \star (b, h) = (abc(g,h) gh)$
then this is an associative product if and only if $c$ satisfies the 
{\em cocycle condition}
$$
c(g,h)c(gh,k) = c(g, hk)c(h,k)
$$
for all $g$, $h$ and $k$ in $G$.  

If we choose a different section $\tilde s$ then we must have
$\tilde s = s d$ for some $d \colon G \to A$.  If $\tilde c$ is the 
cocycle determined by $\tilde s$ then a calculation shows that 
\begin{equation}
    \label{eq:equiv_ext}
    c(g, h) = \tilde c(g, h) d(gh) d(g)^{-1}d(h)^{-1}.
\end{equation}

We have now  essentially shown that all central extensions of $G$ by 
$A$ are determined by cocycles $c$ modulo identifying two 
that satisfy the condition \eqref{eq:equiv_ext}.  Let us recast this 
result in a form that we will see again in this talk. 

Define maps  $d_i \colon G^{p+1}  \to G^{p}$  by 
\begin{equation}
    \label{eq:barcomplex}
d_{i}(g_{1},\ldots,g_{p+1}) = \begin{cases} 
                              (g_{2},\ldots,g_{p+1}), & i = 0, \\ 
                              (g_{1},\ldots,g_{i-1}g_{i},g_{i+1},
                                \ldots,g_{p+1}), & 1\leq i\leq p-1, \\ 
                               (g_{1},\ldots,g_{p}), & i = p. 
                              \end{cases} 
\end{equation}
If $M^p(G;A) = \text{Map}(G^p, A)$ then we define $\delta \colon 
M^p(G; A)\to M^{p-1}(G; A) $ by $\delta(c) = (c\circ d_1) (c \circ 
d_2)^{-1} (c\circ d_3) \dots$. This satisfies $\delta^2 = 1$ and 
defines  a complex
$$
M^0(G;A) \stackrel{\delta}{\to} M^1(G;A)\stackrel{\delta}{\to} M^2(G;A)
\stackrel{\delta}{\to} \dots
$$
The 
cocycle  condition can be rewritten as $\delta(c) = 1$ and the 
condition that two cocycles give rise to the same central extension 
is that $c =\tilde c \delta(d)$.  If we define 
$$
H^p(G; A) =\frac{\text{kernel\ } \delta \colon M^p(G;A) \to M^{p-1}(G;A)}
{\text{image\ } \delta \colon M^{p+1}(G;A) \to M^{p}(G;A)}
$$
then we have shown that central extensions of $G$ by $A$ are classified by 
$H^2(G; A)$. 

\section{Central extensions of Lie groups}
In the case that $G$ is a topological or Lie group it is well-known
that there are interesting central extensions for which no continuous
or differentiable section exists.  For example consider the central 
extension 
$$
\Z_2 \to SU(2) = \text{Spin}(3) \to SO(3) 
$$
of the three dimensional orthogonal group $SO(3)$ by its 
double cover $\text{Spin}(3)$.  Here $SU(2)$ is known to 
be the three sphere but if a section existed then we would have
$SU(2)$ homeomorphic to $\Z_2 \times SO(3)$ and hence disconnected.

From now on we will concentrate on the case when $A = \C^\times$. 
Then $\hat G \to G$ can be thought of as a $\C^\times $ 
principal bundle and a section will only exist if this bundle
is trivial.   The structure of the central extension 
as a $\C^\times$ bundle is important in what follows so we digress
to discuss them in more detail.

\subsection{$\cstar$ bundles}

Let $P \to X$ be a $\cstar$ bundle over a manifold $X$. We  denote
the fibre of $P$ over $x \in X$ by $P_x$. 
Recall \cite{Bry} that if $P$ is a 
 $\cstar$ bundle over a manifold $X$ we can define the dual 
bundle $P^*$ as the same space $P$ 
but with the action $p^* g = (pg^{-1})^*$ and, that if $Q$ is another
such bundle, we can define the product bundle $P\otimes Q$ by 
$(P\otimes Q)_x = (P_x \times Q_x)/\cstar $ where $\cstar$ acts
by $(p,q)w = (pw, qw^{-1})$. We denote an element of $P\otimes Q$ by 
$p \otimes q$ with the understanding that $(pw) \otimes q = p \otimes 
(qw) = (p\otimes q)w$ for $w \in \cstar$. It is straightforward to check that 
$P\otimes P^*$ is canonically trivialised by the section  $x \mapsto p 
\otimes p^*$ where $p$ is any point in $P_x.$

If $P$ and $Q$ are $\cstar$ bundles on $X$ with connections 
$\mu_P$ and $\mu_Q$ then $P \otimes Q$ has an induced
connection we denote by $\mu_P \otimes \mu_Q$. The curvature
of this connection is $R_P + R_Q$ where $R_P$ and $R_Q$ are the 
curvatures of $\mu_P$ and $\mu_Q$ respectively. 
The bundle $P^*$ has an induced connection whose curvature
 is $-R_P$.
 
Recall the maps $d_i \colon G^p \to G^{p-1}$ defined by 
\eqref{eq:barcomplex}.
If $P \to G^p$ is a $\cstar$ bundle then we can define a $\cstar$ 
bundle over $G^{p+1}$ denoted  $\delta(P)$ by 
$$
\delta(P) = \pi_1^{-1}(P) \otimes \pi_2^{-1}(P)^* \otimes \pi_3^{-1}(P)
\otimes \dots.
$$
If $s$ is a section of $P$ then it defines $\delta(s)$ a section of 
$\delta(P)$ and if $\mu$ is a connection on $P$ with curvature
$R$ it defines  a connection on $\delta(P)$ which 
we denote by $\delta(\mu)$.  To define the curvature of $\delta(\mu)$ 
let us  denote by $\Omega^q(G^p)$ the space of all differentiable 
$q$ forms  on $G^p$.  Then we define a map 
\begin{equation}
    \label{delta_forms}
\delta \colon \Omega^q(G^p) \to \Omega^q(G^{p+1}) 
\end{equation}
by $\delta = \sum_{i=0}^p d_i^*$, the alternating sum of pull-backs
by the various maps $d_i \colon G^{p+1} \to G^p$. Then the curvature of 
$\delta(\mu)$ is $\delta(R)$. 
If we consider $\delta(\delta(P))$ it is a product of factors 
and every factor occurs with its dual so $\delta(\delta(P))$ is canonically 
trivial.  If $s $ is a section of $P$ then under this identification 
$\delta\delta(s) = 1$ and moreover if $\mu$ is a connection on $P$ 
then $\delta\delta(\mu) $ is the flat connection on 
$\delta\delta(P)$ with respect to $\delta(\delta(s))$.

\section{Central extensions}
Let $\G$ be a Lie group and
consider a central extension
$$
\cstar\to \hat{\G}\stackrel{\pi}{\to} \G.
$$ 
Following Brylinski and McLaughlin \cite{BryMac}
we think of this as a $\cstar$ bundle $\hat\G \to \G$
with a product $M \colon \hat \G \times \hat  \G
\to \hat\G$ covering the product $m = d_1 \colon G \times G \to G$. 

Because this is a central extension we must have that 
$M(pz, qw) = M(p,q)zw$ for any $p, q \in P$ and $z, w \in 
\cstar$.  This means we have a  section $s$ 
of $\delta(P)$ given by 
$$
s(g, h)  = p \otimes M(p, q) \otimes q
$$
for 
any $p \in P_g$ and $q \in P_h$.  This is well-defined as $pw \otimes 
M(pw, qz) \otimes qz = pw \otimes 
M(p, q)(wz)^{-1} \otimes qz = p \otimes  M(p, q) \otimes q$. 
Conversely any such section gives rise to an $M$. 

Of course we need an associative product and it can be
shown that $M$ being associative is equivalent to $\delta(s) = 1$. 
To actually make $\hat \G$ into a group we need more
than multiplication we need an identity $\hat e \in \hat \G$ and
an inverse map.  It is straightforward to check that if $e \in \G$
is the identity then, because $M \colon \hat \G_e \times \hat \G_e
\to \hat\G_e$,  there is a unique $\hat e \in \hat \G_e$ such that 
$M(\hat e, \hat e) = \hat e$.  It is also straightforward to 
deduce the existence of a unique inverse. 

Hence we have the result from \cite{BryMac} that 
a central extension of $\G$ is a $\cstar$ bundle $P \to G$ 
together with a section $s $ of $\delta(P) \to G\times G$
such that $\delta(s) = 1$.  In \cite{BryMac} this is phrased
in terms of simplicial line bundles.  Note that this is a kind of 
cohomology result analogous to that in the first section. We have an 
object (in this case a $\Cstar$ bundle) and $\delta$ of the 
object is `zero' i.e. trivial as a $\Cstar$ bundle. 

For our purposes we need to phrase this result in terms of 
differential forms. 
We call a connection for $\hat \G \to \G$, thought of as a 
$\cstar$ bundle, a connection for the central extension.  An 
isomorphism of central extensions with connection is an isomorphism
of bundles with connection which is a group isomorphism on the total
space $\hat\G$.  Denote by 
$C(\G)$ the set of all isomorphism classes of 
central extensions of $\G$ with connection.

Let  $\mu \in \Omega^1(\hat \G)$  be a connection on the
bundle  $\hat{\G} \to \G$ and consider the tensor product connection 
$\delta(\mu)$. Let
 $\alpha = s^*(\delta(\mu))$. 
We then have that 
\begin{align*}
    \delta(\alpha) &= (\delta(s)^*) (\delta(\mu)) \\
&= (1)^*(\delta^2(\mu))\\
&= 0
\end{align*}
as $\delta^2(\mu) $ 
is the flat connection on $\delta^2(P)$.   Also $d\alpha = 
s^*(d\delta(\mu)) = \delta(R)$. 

Let $\Gamma(\G)$ denote the set of all pairs $(\alpha, R)$ where 
$R$ is a closed, $2\pi i $ integral, two form on $\G$ and $\alpha$ is a 
one-form on $\G \times \G$ with $\delta(R) = d\alpha$ and 
$\delta(\alpha) = 0$. 

We have constructed  
a map $C(\G) \to \Gamma(\G)$.
In the next section we construct an inverse to this 
map by showing how to define a central extension from a pair 
$(\alpha, R)$.  For now notice that isomorphic central 
extensions with connection clearly give rise to the 
same $(\alpha, R)$ and that if we vary the connection, which 
is only possible by adding on the pull-back of a 
one-form $\eta$ from $\G$, then we change $(\alpha, R)$
to $(\alpha + \delta(\eta), R + d\eta)$.

\subsection{Constructing the central extension}
Recall that given $R$ we can find a principal $\cstar$ bundle $P \to 
\G$  with connection $\mu$ and curvature $R$ which is 
unique up to isomorphism. 
It is a standard result in the theory of bundles that 
if $P \to X$ is a bundle with connection $\mu$ which 
is flat and $\pi_1(X) = 0$ then $P$ has a section $s \colon 
X \to P$ such that $s^*(\mu) = 0$. Such a section is not 
unique of course it can be multiplied by a (constant) element of 
$\cstar$.  As our interest is in the loop group $G$ which satisfies $\pi_1(G)= 0$
we shall assume, from now on, that $\pi_1(G) = 0$. 
Consider now our pair $(R, \alpha)$ and the bundle $P$. 
As $\delta(R) = d\alpha$ we have that the connection $\delta(w) - \pi^*(\alpha)$
on $\delta(P) \to G\times G$ is flat and hence (as $\pi_1(G \times 
G) = 0$) we can find a section $s$  such that 
$s^*(\delta(w)) = \alpha$.

The section $s$ defines a multiplication by 
$$
s(p, q) = p \otimes M(p, q)^* \otimes q.
$$
Consider now $\delta(s)$ this satisfies 
$\delta(s)^*(\delta(\delta(w))) = \delta(s^*(\delta(w)) = \delta(\alpha) = 
0$.
 On the other hand the canonical section $1$ of $\delta(\delta(P))$
also satisfies this so they differ by a constant element of the group. 
This means that there is a $w \in \cstar$ such that for any $p$, $q$ and $r$ we must have 
$$
M(M(p, q), r) = w M(p, M(q, r)).
$$
Choose $p \in \hat\G_e$ where $e$ is the identity in $\G$. Then
$M(p, p ) \in \hat\G_e$ and hence $M(p, p) = pz$ for some $z \in 
\cstar$.  Now let $p=q=r$ and it is clear that we must have $w=1$. 

So from $(\alpha, R)$ we have constructed $P$ and a section $s$  of 
$\delta(P)$ with $\delta(s) = 1$.  However $s$ is not unique
but this is not a problem. 
If we change $s$  to $s' = sz$ for 
some constant $z \in \cstar$ then we have changed $M$
to $M' = M z$. As $\cstar$ is central multiplying
by $z$ is an isomorphism of central extensions with 
connection. So the ambiguity in $s$  does not change the 
isomorphism class of the central extension with connection. Hence
we have constructed a map 
$$
\Gamma(\G) \to C(\G)
$$
as required. That it is the inverse of the earlier
map follows from the definition of $\alpha$ as $s^*(\delta(\mu))$
and the fact that the connection on $P$ is chosen so its 
curvature is $R$.

\section{Conclusion: Loop groups}

In the case where $\G = L (K)$ there is a 
well known expression for the curvature $R$ of a 
left invariant connection on $\hat{L(K)}$ 
--- see \cite{PreSeg}.  We can also write down 
a 1-form $\alpha$ on $L(K)\times L(K)$ 
such that $\d(R) = d\alpha$ and $\d(\alpha) = 0$.  We have:
\begin{align*} 
 R(g)(gX,gY) &= \frac{1}{4\pi^{2}}\int_{S^1} \<
X , \partial_\theta Y \> d\theta    \\ 
 \alpha(g_{1},g_{2})(g_{1}X_{1},g_{2}X_{2}) &= \frac{1}{4\pi^{2}} 
\int_{S^{1}} \< X_{1} ,(\partial_\theta g_{2}){\partial_\theta} g_{2}^{-1} \> d\theta.  
\end{align*} 
Here $\<\ ,\ \>$ is the Killing form on $\mathfrak{k}$ normalised
so the longest root has length squared equal to $2$ and 
$\partial_\theta$ denotes differentiation with respect to $\theta \in 
S^1$.
Note that $R$ is left invariant and that $\alpha$ 
is left invariant in the first factor of  $\G\times \G$.  It can be 
shown that these are the $R$ and $\alpha$ arising in \cite{Mur}.

In \cite{MurSte} we apply the methods of this talk to give an 
explicit construction of the `string class' of a loop group 
bundle and relate it to earlier work of Murray on calorons.

\end{document}